\documentclass{article}\begin{document}
\centerline{\bf\large On Algebraic Functions}\vskip .2in

\centerline{N.D. Bagis}

\centerline{Stenimahou 5 Edessa}
\centerline{Pella 58200, Greece}
\centerline{bagkis@hotmail.com}

\[
\]

\centerline{\bf Abstract}

\begin{quote}
In this article we consider functions with Moebius-periodic rational coefficients. These functions under some conditions take algebraic values and can be recovered by theta functions and the Dedekind eta function. Special cases are the elliptic singular moduli, the Rogers-Ramanujan continued fraction, Eisenstein series and functions associated with Jacobi symbol coefficients.         
\end{quote}


\section{Known Results on Algebraic Functions}

The elliptic singular moduli $k_r$ is the solution $x$ of the equation
\begin{equation}
\frac{{}_2F_{1}\left(\frac{1}{2},\frac{1}{2};1;1-x^2\right)}{{}_2F_{1}\left(\frac{1}{2},\frac{1}{2};1;x^2\right)}
=\sqrt{r}
\end{equation}
where 
\begin{equation}
{}_2F_{1}\left(\frac{1}{2},\frac{1}{2};1;x^2\right)=\sum^{\infty}_{n=0}\frac{\left(\frac{1}{2}\right)^2_n}{(n!)^2} x^{2n}=\frac{2}{\pi}K(x)=\frac{2}{\pi}\int^{\pi/2}_{0}\frac{d\phi}{\sqrt{1-x^2\sin^2(\phi)}}
\end{equation}
The 5th degree modular equation which connects $k_{25r}$ and $k_r$ is (see [13]):
\begin{equation}
k_rk_{25r}+k'_rk'_{25r}+2^{5/3} (k_rk_{25r}k'_rk'_{25r})^{1/3}=1
\end{equation} 
The problem of solving (3) and find $k_{25r}$ reduces to that solving the depressed equation after named by Hermite (see [3]):
\begin{equation}
u^6-v^6+5u^2v^2(u^2-v^2)+4uv(1-u^4v^4)=0
\end{equation}
where $u=k^{1/4}_{r}$ and $v=k^{1/4}_{25r}$.\\
The function $k_r$ is also connected to theta functions from the relations 
\begin{equation}
k_r=\frac{\theta^2_2(q)}{\theta^2_3(q)}, \textrm{ where } \theta_2(q)=\theta_2=\sum^{\infty}_{n=-\infty}q^{(n+1/2)^2} \textrm { and } \theta_3(q)=\theta_3=\sum^{\infty}_{n=-\infty}q^{n^2}
\end{equation}
$q=e^{-\pi\sqrt{r}}$.\\  
Hence a closed form solution of the depressed equation is
\begin{equation}
k_{25r}=\frac{\theta^2_2(q^5)}{\theta^2_3(q^5)}
\end{equation} 
But this is not satisfactory.\\ For example in the case of $\pi$ formulas of Ramanujan (see [12] and related references), one has to know from the exact value of $k_r$ the exact value of $k_{25r}$ in radicals. (Here we mention the concept that when $r$ is positive rational then the value of $k_r$ is algebraic number). Another example is the Rogers-Ramanujan continued fraction (RRCF) which is
\begin{equation}
R(q)=\frac{q^{1/5}}{1+}\frac{q}{1+}\frac{q^2}{1+}\frac{q^3}{1+}\ldots
\end{equation}
(see [4],[5],[6],[8],[9],[10],[13],[14],[16],[21]),   
the value of which depends form the depressed equation.\\ If we know the value of (RRCF) then we can find the value of $j$-invariant from Klein's equation (see [19],[8] and Wolfram pages 'Rogers Ramanujan Continued fraction'):
\begin{equation}
j_r=-\frac{\left(R^{20}-228 R^{15}+494 R^{10}+228 R^5+1\right)^3}{R^5 \left(R^{10}+11 R^5-1\right)^5} \textrm{ , where } R=R(q^2)
\end{equation}
One can also prove that Klein's equation (8) is equivalent to depressed equation (4).\\ Using the 5th degree modular equation of Ramanujan
\begin{equation}
R(q^{1/5})^5=R(q)\frac{1-2R(q)+4R(q)^2-3R(q)^3+R(q)^4}{1+3R(q)+4R(q)^2+2R(q)^3+R(q)^4}
\end{equation} 
and (8) we can find the value of $j_{r/25}$ and hence from the relation 
\begin{equation}
j_r=\frac{256(k^2_r+k'^4_{r})^3}{(k_rk'_r)^4} . 
\end{equation}
$k_{r/25}$. Knowing $k_{r}$ and $k_{r/25}$, we have evaluated $k_{25r}$ (see [7]) and give relations of the form 
\begin{equation}
k_{25r}=\Phi(k_r,k_{r/25})\textrm{ and }k_{25^nr}=\Phi_n(k_r,k_{r/25}), n\in\bf N\rm 
\end{equation}
Hence when we know the value of $R(q)$ in radicals we can find $k_r$ and $k_{25r}$ in radicals and the opposite.\\
\\
Also in [3] and Wikipedia 'Bring Radical' one can see how the depressed equation can used for the extraction of the solution of the general quintic equation
\begin{equation}
ax^5+bx^4+cx^3+dx^2+ex+f=0 
\end{equation}
The above equation can solved exactly with theta functions and in some cases in radicals.\\ The same holds and with the sextic equation (see [8])
\begin{equation}
\frac{b^2}{20a}+bY+aY^2=cY^{5/3} 
\end{equation}
which have solution
\begin{equation}
Y=Y_r=\frac{b}{250a}\left(R(q^2)^{-5}-11-R(q^2)^5\right)\textrm{, }q=e^{-\pi\sqrt{r}}\textrm{, }r>0
\end{equation}
and $r$ can evaluted from the constants using the relation  $j_r=250\frac{c^3}{a^2b}$, in order to generate the solution. 
\\
The Ramanujan-Dedekind eta function is defined as
\begin{equation}
\eta(\tau)=\prod^{\infty}_{n=1}(1-q^n)\textrm{, }q=e^{i\pi\tau}\textrm{, }\tau=\sqrt{-r}
\end{equation}  
The $j$-invariant can be expresed in terms of Ramanujan-Dedekind eta function as
\begin{equation}
j_r=\left[\left(q^{-1/24}\frac{\eta(\tau)}{\eta(2\tau)}\right)^{16}+16\left(q^{1/24}\frac{\eta(2\tau)}{\eta(\tau)}\right)^{8}\right]^3
\end{equation}
The Ramanujan-Dedekind eta function is (see [11],[22])
\begin{equation}
\eta(\tau)^8=\frac{2^{8/3}}{\pi^4}q^{-1/3}(k_r)^{2/3}(k'_r)^{8/3}K(k_r)^4
\end{equation} 
There are many interesting things one can say about algebricity and special functions.\\
In this article we examine Moebius-periodic functions. If the Taylor coefficients of a function are Moebius periodic, then we can evaluate prefectly these functions, by taking numerical values. This can be done using the Program Mathematica and the routine 'Recognize'. By this method we can find values coming from the middle of nowhere. However they still remaining conjectures. Also it is great challenge to find the polynomials and modular equations of these Moebius-periodic functions and united them with a general theory. Many various scientists are working for special functions (mentioned above) such as the singular moduli ($j$-invariant) and the related to them Hilbert polynomials, (RRCF), theta functions, Dedekind-Ramanujan $\eta$ and other similar to them. In [7] is presented a way to evaluate the fifth singular moduli and the Rogers-Ramanujan continued fraction with the function $w_r=\sqrt{k_rk_{25r}}$. This function can replace the classical singular moduli in the case of Rogers-Ramanujan continued fraction and Klein's invariant.\\ We are concern to construct a theory of such functions and characterize them.

\section{The Main Theorem}

We begin by giving a definition and a conjecture which will help us for the proof of the Main Theorem.\\ 
\\
\textbf{Definition 1.}\\Let $a$,$p$ be positive rational numbers with $a<p$ and $q=e^{-\pi\sqrt{r}}$, $r>0$. We call ''agiles'' the quantities 
\begin{equation}
[a,p;q]:=\prod^{\infty}_{n=0}(1-q^{pn+a})(1-q^{pn+p-a})
\end{equation}
\\
The ''agiles'' have the following very interesting conjecture-property.\\
\\
\textbf{Conjecture.}\\ If $q=e^{-\pi\sqrt{r}}$, $r$ is positive rational and $a,b$ are positive rationals with $a<p$, then 
\begin{equation}
[a,p;q]^{*}:=q^{p/12-a/2+a^2/(2p)}[a,p;q]=\textrm{Algebraic Number}
\end{equation}
\\
Assuming the above unproved property we will show the following\\
\\
\textbf{Main Theorem.}\\
Let $f$ be a function analytic in $(-1,1)$. Set by Moebius theorem bellow $X(n)$ to be
\begin{equation}
X(n)=\frac{1}{n}\sum_{d|n}\frac{f^{(d)}(0)}{\Gamma(d)}\mu\left(\frac{n}{d}\right)
\end{equation}  
If $X(n)$ is $T$ periodic sequence, rational valued and catoptric in the every period-interval i.e. if for every $n\in\bf N\rm$ is $a_k=X(k+nT)=X(k)$ with $a_T=0$ we have $a_1=a_{T-1}, a_2=a_{T-2},\ldots,a_{(T-1)/2}=a_{(T+1)/2}$, then exist rational a number $A$, such that
\begin{equation}   
q^Ae^{-f(q)}=\textrm{Algebraic Number} , 
\end{equation} 
in the points $q=e^{-\pi\sqrt{r}}$, $r$ positive rational.\\
The number $A$ is given from  
\begin{equation}
A=\sum^{\left[\frac{T-1}{2}\right]}_{j=1}\left(-\frac{j}{2}+\frac{j^2}{2T}+\frac{T}{12}\right)X(j)
\end{equation}
\\

For to prove the Main Theorem we will use the next known (see [2]) Moebius inversion Theorem\\
\\
\textbf{Theorem.} (Moebius inversion Theorem)\\
If $f(n)$ and $g(n)$ are arbitrary arithmetic functions, then
\begin{equation}
\sum_{d|n}f(d)=g(n)\Leftrightarrow f(n)=\sum_{d|n}g(d)\mu\left(\frac{n}{d}\right)
\end{equation}
\\
The Moebius $\mu$ function is defined as $\mu(n)=0$ if $n$ is not square free, and $(-1)^r$ if $n$ have $r$ distinct primes.\\
\\ 
Hence some values are $\mu(1)=1$, $\mu(3)=-1$, $\mu(15)=1$, $\mu(12)=0$, etc.\\
For to prove the Main Theorem we will use the next\\
\\
\textbf{Lemma 1.}\\
If $|x|<1$
\begin{equation}
\log\left(\prod^{\infty}_{n=1}\left(1-x^n\right)^{X(n)}\right)=-\sum^{\infty}_{n=1}\frac{x^n}{n}\sum_{d|n}X(d)d
\end{equation} 
\\
\textbf{Proof.}\\ It is $|x|<1$, hence
$$
\log\left(\prod^{\infty}_{n=1}\left(1-x^n\right)^{X(n)}\right)=\sum^{\infty}_{n=1}X(n)\log\left(1-x^n\right)=
$$
$$
=-\sum^{\infty}_{n=1}X(n)\sum^{\infty}_{m=1}\frac{x^{mn}}{m}=-\sum^{\infty}_{n,m=1}\frac{x^{nm}}{nm}X(m)m=-\sum^{\infty}_{n=1}\frac{x^n}{n}\sum_{d|n}X(d)d .
$$
\\
\textbf{Proof of Main Theorem.}\\
From Taylor expansion theorem we have
$$
e^{-f(x)}=\exp\left(-\sum^{\infty}_{n=1}\frac{f^{(n)}(0)}{n!}x^n\right)
$$
From the Moebius inversion theorem exists $X(n)$ such that 
$$
X(n)=\frac{1}{n}\sum_{d|n}\frac{f^{(d)}(0)}{\Gamma(d)}\mu\left(\frac{n}{d}\right)
$$  
or equivalent
$$
\frac{f^{(n)}(0)}{n!}=\frac{1}{n}\sum_{d|n}X(d)d
$$
Hence from Lemma 1
$$
e^{-f(x)}=\exp\left(-\sum^{\infty}_{n=1}\frac{x^n}{n}\sum_{d|n}X(d)d\right)=\prod^{\infty}_{n=1}\left(1-x^n\right)^{X(n)}
$$
and consequently because of the periodicity and the catoptric property of $X(n)$, we get
\begin{equation}
e^{-f(x)}=\prod^{\left[\frac{T-1}{2}\right]}_{j=1}[j,T;x]^{X(j)}
\end{equation}
which is a finite product of ''agiles'' and from the Conjecture exist $A$ rational such that (21) hold, provided that $x=q^{-\pi\sqrt{r}}$ and $r$ positive rational.\\    
\\
\textbf{Examples.}\\

\textbf{1)} For $X(n)=\left(\frac{n}{G}\right)$, where $G=2^mg_1^{m_1}g_2^{m_2}\ldots g_s^{m_s}$, $m,m_1,\dots,m_s$ not negative integers, $m\neq1$ and $g_1<g_2<\ldots<g_s$ primes of the form $1\textrm{mod}4$, then exist $A$ rational such that
$$
q^A\prod^{\infty}_{n=1}(1-q^n)^{\left(\frac{n}{G}\right)}=q^A\prod^{\left[\frac{G-1}{2}\right]}_{j=1}[j,G,q]^{X(j)}=\textrm{Algebraic}
$$
when $q=e^{-\pi\sqrt{r}}$, $r$ positive rational.\\
A special case is $G=5$ which gives the Rogers-Ramanujan continued fraction. More precisely is $X=\{1,-1,-1,1,0,...\}$ and evaluations can given.
$$
q^{1/5}\prod^{\infty}_{n=1}(1-q^n)^{\left(\frac{n}{5}\right)}=R(q)=q^{1/5}\frac{1}{1+}\frac{q}{1+}\frac{q^2}{1+}\ldots
$$

\textbf{2)} If $X=\{1,1,0,1,1,0,\ldots\}$, then $T=3$ and
$A=-\frac{1}{12}$. Hence we get that if $q=e^{-\pi}$, then  
$$
q^{-1/12}e^{-f(q)}=\sqrt[12]{81 \left(885+511 \sqrt{3}-3 \sqrt{174033+100478 \sqrt{3}}\right)}
$$

\textbf{3)} If $X=\{1,1,1,1,0,1,1,1,1,0,\ldots\}$, then $T=5$ and $A=\frac{-1}{6}$ we get that\\ 
i) If $q=e^{-\pi\sqrt{2}}$, then $q^{-1/6}e^{-f(q)}$, is root of
$$
3125+250 v^6-20 v^{10}+v^{12}=0
$$
We can solve the above equation observing that 
is of the form (13): 
$$
3125+250Y_r^6+Y_r^{12}=j_r^{1/3}Y_r^{10},\eqno{(eq)}
$$ 
where $j_r$ is the $j$-invariant. Hence 
$$
q^{-1/6}e^{-f(q)}=\sqrt[6]{Y_{1/2}}
$$
then see [21]:
$$
R(e^{-2\pi\sqrt{2}})=\frac{\sqrt{5(g+1)+2g\sqrt{5}}-\sqrt{5g}-1}{2}
$$
where
$$
(g^3-g^2)/(g+1)=(\sqrt{5}+1)/2
$$
One can use the duplication formula of RRCF (see [16]) to find $R(e^{-\pi\sqrt{2}})$ in radicals and hence the value of $Y_{1/2}$ in radicals.\\
ii) If $q=e^{-2\pi}$, then $$q^{-1/6}e^{-f(q)}=\sqrt[6]{Y_{1}}=\sqrt{\frac{5}{2}+\frac{5 \sqrt{5}}{2}}$$
...etc\\
If $q=e^{-\pi\sqrt{r}}$ 
$$
q^{-1/6}e^{-f(q)}=\sqrt[6]{Y_{r/4}}
$$

\section{The Representation of $e^{-f(q)}$}

We give in Theorem 1 bellow the representation of a Moebius periodic function $f$ in terms of known functions.\\    
For $|q|<1$ the Jacobi theta functions are 
\begin{equation}
\vartheta(a,b;q):=\sum^{\infty}_{n=-\infty}(-1)^nq^{an^2+bn}
\end{equation}
One can evaluate the agiles by the theta function
\begin{equation}
M(c,q):=\sum^{\infty}_{n=0}c^nq^{n(n+1)/2}=\frac{1}{1+}\frac{-cq}{1+}\frac{-c(q-q^2)}{1+}\frac{-cq^3}{1+}\frac{-c(q^4-q^2)}{1+}\ldots
\end{equation}
then from [11] a way to express the agiles is
\begin{equation}
[a,p;q]=\frac{M(-q^{-a},q^p)-q^aM(-q^a,q^p)}{\eta(\tau p)}
\end{equation}
however we shall use the general theta functions evaluation see relation (32) bellow for more concentrated forms. The reader can change from one form to the other.\\ 
\\
A first result in the agiles also given in [11] was the evaluation of the duplication formula 
\begin{equation}
\frac{[a,p;q^2]^{*}}{[a,p;q]^{*}}=\tau^{*}(a,p;q),
\end{equation}
for which if $a,b$ are positive reals and $n$ integer, then
\begin{equation}
\tau^{*}(a,p;q)=\tau^{*}(np\pm a,p;q)
\end{equation}
\\
\textbf{Theorem 1.}\\
For every $f$ analytic in $(-1,1)$ with 
$$X(n)=\frac{1}{n}\sum_{d|n}\frac{f^{(d)}(0)}{\Gamma(d)}\mu\left(\frac{n}{d}\right)$$
periodic-symmetric with period $T$ (Moebius periodic) and real-valued, then hold
\begin{equation}
e^{-f(q)}=\eta(T\tau)^{-\sum^{\left[\frac{T-1}{2}\right]}_{j=1}X(j)}\prod^{\left[\frac{T-1}{2}\right]}_{j=1}\vartheta\left(\frac{T}{2},\frac{T-2j}{2};q\right)^{X(j)}
\end{equation}
for every $|q|<1$. In case that $X(j)$ are rational, then $q^Ae^{-f(q)}$ is algebraic when $q=e^{-\pi\sqrt{r}}$, $r$ positive rational.\\ 
\\
\textbf{Proof.}\\
Use the expansion found in [11]:
\begin{equation}
[a,p;q]=\frac{1}{\eta(p\tau)}\sum^{\infty}_{n=-\infty}(-1)^nq^{pn^2/2+(p-2a)n/2}=\frac{1}{\eta(p\tau)}\vartheta\left(\frac{p}{2},\frac{p-2a}{2};q\right)
\end{equation}
along with relation (25).\\ 
\\
\textbf{Theorem 2.}\\
If $X(n)$ is real $T$-periodic and catoptric, then
\begin{equation}
\sum^{\infty}_{n=1}\frac{nX(n)q^n}{1-q^n}=-q\frac{d}{dq}\log\left(\eta(T\tau)^{-\sum^{\left[\frac{T-1}{2}\right]}_{j=1}X(j)}\prod^{\left[\frac{T-1}{2}\right]}_{j=1}\vartheta\left(\frac{T}{2},\frac{T-2j}{2};q\right)^{X(j)}\right)
\end{equation}
\textbf{Proof.}\\
If 
$$
X(n)=\frac{1}{n}\sum_{d|n}\frac{f^{(d)}(0)}{\Gamma(d)}\mu(n/d)$$
then it holds
$$
\int\frac{1}{q}\sum^{\infty}_{n=1}\frac{nX(n)q^n}{1-q^n}dq=f(q)
$$
From Theorem 1 we get the result.\\
\\
\textbf{Remark.}\\
If $R(a,b,p;q)=q^C\frac{[a,p;q]}{[b,p;q]}$ denotes a Ramanujan quantity (see[5]), then we have the next closed form evaluation, with theta functions
\begin{equation}
\frac{R'\left(X_p,q\right)}{R\left(X_p,q\right)}=\frac{C}{q}+\frac{d}{dq}\log\left(\prod^{\left[\frac{p-1}{2}\right]}_{j=1}\vartheta\left(\frac{p}{2},\frac{p-2j}{2},q\right)^{X_p(j)}\right)
\end{equation}
by using (33).

\section{Examples and Applications}

\textbf{4)} The Jacobi symbol $\left(\frac{n}{5}\right)$ is 5 periodic and symmetric, hence
\begin{equation}
\sum^{\infty}_{n=1}\left(\frac{n}{5}\right)\frac{nq^n}{1-q^n}=-q\frac{d}{dq}\log\left(\frac{\vartheta(5/2,3/2;q)}{\vartheta(5/2,1/2;q)}\right)=-q\frac{d}{dq}\log(q^{-1/5}R(q)) .
\end{equation}
\\
From Example 3 and the above relation (35), one can prove that 
\begin{equation}
\sum^{\infty}_{n=1}\frac{nq^n}{1-q^n}-\sum^{\infty}_{n=1}\frac{5nq^{5n}}{1-q^{5n}}=-q\frac{d}{dq}\log\left(\frac{\vartheta(5/2,1/2;q)\vartheta(5/2,3/2;q)}{\eta\left(5\tau \right)^{2}}\right)
\end{equation}
and
\begin{equation}
\frac{1}{6}+\sum^{\infty}_{n=1}\frac{nq^n}{1-q^n}-5\sum^{\infty}_{n=1}\frac{nq^{5n}}{1-q^{5n}}=-\frac{q}{6}\frac{d}{dq}\log\left(Y\left(\sqrt{q}\right)\right)
\end{equation}
In view of [5] relation (92) and the expansion of $L_1(q)$ in the same paper, we get 
$$
-\frac{q}{6}\frac{d}{dq}\log\left(Y\left(q^{1/2}\right)\right)=\frac{-q^{1/2}}{12}\frac{Y'\left(q^{1/2}\right)}{Y\left(q^{1/2}\right)}
=-\frac{1}{6}-\frac{K[r]^2}{6 \pi^2}+\frac{a(r) K[r]^2}{\pi^2 \sqrt{r}}+\frac{5 K[25 r]^2}{6 \pi^2}-
$$
\begin{equation}
-\frac{a(25 r) K[25 r]^2}{\pi^2 \sqrt{r}}-\frac{K[r]^2 k_r^2}{6 \pi^2}+\frac{5 K[25 r]^2 k_{25 r}^2}{6\pi^2}
\end{equation}
where $a(r)$ is the elliptic alpha function (see [17]) and $K[r]=K(k_r)$ is the complete elliptic integral of the first kind at singular values. Hence for certain $r$ we can find special values of  $Y'\left(e^{-\pi\sqrt{r}}\right)$.\\
   
From (36) and (37) we get the following evaluation for the theta function
\begin{equation}
\theta=\frac{\vartheta(5,1;q)^6\vartheta(5,3;q)^6}{q^2\eta\left(10\tau \right)^{12}}=R(q^2)^{-5}-11-R(q^2)^5
\end{equation}
and in view of [9] we get the following, similar to inverse elliptic nome theorem\\
\\
\textbf{Theorem 3.}
\begin{equation}
\frac{-1}{5}\int^{\theta}_{+\infty}\frac{dt}{t^{1/6}\sqrt{125+22t+t^2}}=\frac{1}{5\sqrt[3]{4}}B(k_{4r},1/6,2/3)
\end{equation}
and
\begin{equation}
\theta=H\left(k_{4r}\right)
\end{equation}
where
\begin{equation}
\frac{-1}{5}\int^{G(x)}_{+\infty}\frac{dt}{t^{1/6}\sqrt{125+22t+t^2}}=x\textrm{ and }H(x)=G\left(\frac{B\left(x,\frac{1}{6},\frac{2}{3}\right)}{5\sqrt[3]{4}}\right)
\end{equation}
Also 
\begin{equation}
\frac{d}{dr}B(k_r^2,1/6,2/3)=-\frac{\pi}{2}\sqrt[3]{4}\frac{q^{1/6}\eta\left(\tau\right)^4}{\sqrt{r}}
\end{equation}
and $\theta^{1/6}$ is root of $(eq)$.\\
\\

For example with $r=1/5$, then $\theta=5\sqrt{5}$ and 
$$
k_{4/5}=\frac{2-\sqrt{2-4\sqrt{-2+\sqrt{5}}}}{2+\sqrt{2-4\sqrt{-2+\sqrt{5}}}}
$$
$$
5\sqrt{5}=H\left(\frac{2-\sqrt{2-4\sqrt{-2+\sqrt{5}}}}{2+\sqrt{2-4\sqrt{-2+\sqrt{5}}}}\right) .
$$
\\

Continuing we have, if $X(n)=\left(\frac{n}{G_0}\right)$ have period $p$ and $G_0=2^{m_0}p_1^{m_1}p_2^{m_2}\ldots p_s^{m_s}$, with $p_j$-primes of the form $1(mod4)$, $m_s,s,j=0,1,2,\ldots$ and $m_0\neq 1$, then
\begin{equation}
\sum^{\infty}_{n=1}\left(\frac{n}{G_0}\right)\frac{nq^n}{1-q^n}=-q\frac{d}{dq}\log\left(\prod^{\left[\frac{G_0-1}{2}\right]}_{j=1}\vartheta\left(\frac{G_0}{2},\frac{G_0-2j}{2};q\right)^{\left(\frac{j}{G_0}\right)}\right)
\end{equation}
Also\\ 
\\
\textbf{Conjecture 2.}\\
If $g$ is perfect square and $p_1<p_2<\ldots<p_{\lambda}$ are the distinct primes in the factorization of $g$, then 
\begin{equation}
\prod^{\infty}_{n=1}(1-q^n)^{\left(\frac{n}{g}\right)}
=\eta(\tau)\prod^{\lambda}_{i=1}\eta(p_i\tau)^{-1}\prod_{i<j}\eta(p_ip_j\tau)^1\prod_{i<j<k}\eta(p_ip_jp_k\tau)^{-1}\ldots
\end{equation}
and
$$
\prod^{\infty}_{n=1}(1-q^n)^{-\left(\frac{n}{g}\right)}\frac{d}{dq}\prod^{\infty}_{n=1}(1-q^n)^{\left(\frac{n}{g}\right)}=\sum^{\infty}_{n=1}\left(\frac{n}{g}\right)\frac{n q^n}{1-q^n}=
$$
$$
=-q\frac{d}{dq}\log\left[\eta(g\tau)^{-\sum^{\left[\frac{g-1}{2}\right]}_{j=1}\left(\frac{j}{g}\right)}\prod^{\left[\frac{g-1}{2}\right]}_{j=1}\vartheta\left(\frac{g}{2},\frac{g-2j}{2};q\right)^{\left(\frac{j}{g}\right)}\right]=
$$
$$
=-q^{-1}\left[L(q)-\sum^{\lambda}_{i=1}p_iL(q^{p_i})+\sum_{i<j}p_ip_jL(q^{p_ip_j})-\ldots\right]
$$
where $\eta(p\tau)$ and $L(q^p)$ can evaluated explicitly from (17) and 
\begin{equation}
1-24\sum^{\infty}_{n=1}\frac{nq^n}{1-q^n}
=\frac{6}{\pi\sqrt{r}}+4\frac{K^2[r]\left(-6 \alpha(r)+\sqrt{r}\left(1+k^2_r\right)\right)}{\pi^2 \sqrt{r}},
\end{equation}
$\alpha(r)$ is the elliptic alpha function, $k_r$ the elliptic singular moduli.\\
\\
We also have\\
\\
\textbf{Theorem 4.}\\
If $p$ is prime then
$$
\frac{\pi^{2}\sqrt{r}}{4K[r]^2}\left[-1+p-24q\frac{d}{dq}\log\left(\eta(p\tau)^{-\frac{p-1}{2}}\prod^{\left[\frac{p-1}{2}\right]}_{j=1}\vartheta\left(\frac{p}{2},\frac{p-2j}{2};q\right)\right)\right]=
$$
$$
=6 \alpha(r)-\sqrt{r}\left(1+k^2_r\right)
+m_{p^2r}^2\left(-6 \alpha(p^2r)+p\sqrt{r}\left(1+k^2_{p^2r}\right)\right).
$$
where $m_{n^2r}=\frac{K(k_{n^2r})}{K(k_r)}$ is the multiplier  and is algebraic valued function when $r,n$ are in $\bf Q^{*}_{+}\rm$, (see [13],[17]).

\section{The corespondence of theta functions and singular modulus}

In [11] we have shown that if $m$ is integer and $q=e^{-\pi\sqrt{r}}$, $r>0$ then
\begin{equation}
\sum^{\infty}_{n=-\infty}q^{n^2+2mn}=q^{-m^2}\sqrt{\frac{2K[r]}{\pi}}
\end{equation} 
which is a classical result (see [3]). Also we have shown that
\begin{equation}
\sum^{\infty}_{n=-\infty}q^{n^2+(2m+1)n}=2^{5/6}q^{-(2m+1)^2/4}\frac{(k_{11}k_{12}k_{21})^{1/6}}{k_{22}^{1/3}}\sqrt{\frac{K[r]}{\pi}}
\end{equation} 
where $k_{11}=k_r$, $k_{12}=\sqrt{1-k_{11}^2}$, $k_{21}=\frac{2-k_{11}^2-2k_{12}}{k_{11}^2}$, $k_{22}=\sqrt{1-k_{12}^2}$.\\   
The above relations (47) and (48) give us evaluations of all
\begin{equation}
\sum^{\infty}_{n=-\infty}q^{n^2+mn}
\end{equation}   
with $m$ integer.\\
We define $k_i(x)$ the inverse function of the singular modulus $k_x$. Then it must holds from relation (1)
\begin{equation}
k_i(x)=\left(\frac{K\left(\sqrt{1-x^2}\right)}{K(x)}\right)^2
\end{equation}
From all the above algebraic properties of the theta functions of the previous paragraphs one might think what will happen if we replace $r$ with $k_i\left(\frac{m}{n}\right)$, where $m,n$ positive integers? Is there a nome that says that theta functions are algebraic functions of the singular modulus? The answer is ''yes'' for a given theta function 
$$
[a,p;q]^{*}=\frac{q^{\frac{p}{12}-\frac{a}{2}+\frac{a^2}{2p}}}{\eta(p\tau)}\sum^{\infty}_{n=-\infty}(-1)^nq^{pn^2/2+(p-2a)n/2}=
$$
\begin{equation}
=\frac{q^{\frac{p}{12}-\frac{a}{2}+\frac{a^2}{2p}}}{\eta(p\tau)}\vartheta\left(\frac{p}{2},\frac{p-2a}{2};q\right),
\end{equation}
there exists a unique algebraic function $Q_{\{a,p\}}(x)$ such that    
\begin{equation}
[a,p;e^{-\pi\sqrt{k_i(x)}}]=Q_{\{a,p\}}(x)
\end{equation}
and $Q_{\{a,p\}}(x)$ has the property that is root of a polynomial of a fixed degree depending only on $a$ and $p$.\\
For example the theta function $[1,4,q]$ has  
\begin{equation}
Q_{\{1,4\}}(x)=\sqrt[12]{\frac{4(1-x^2)}{x}}
\end{equation}
Also the function $[1/2,4,q]$ has 
\begin{equation}
Q_{\{1/2,4\}}(x)=\sqrt[48]{\frac{4(1-x)^4(2+x-2\sqrt{1+x})^{12}}{x^{13}(1+x)^2}}
\end{equation}
Hence one may lead to the evaluations   
\begin{equation}
[1,1/2,q]^{*}=\sqrt[12]{\frac{4(1-k_r^2)}{k_r}}
\end{equation} 
and
\begin{equation}
\left[1/2,4,q\right]^{*}=\sqrt[48]{\frac{4\left(1-k_r\right)^4\left(2+k_r-2\sqrt{1+k_r}\right)^{12}}{k_r^{13}\left(1+k_r\right)^2}}
\end{equation}
for every $r>0$.
Other numerical results can also show us that this is conjecturaly true. The function $Q_{\{a,p\}}(x)$ is always algebraic and its degree is the same for all the values of $q_x=e^{-\pi\sqrt{k_i(x)}}$, $x$ positive rational and
\begin{equation}
[a,p;q]^{*}=Q_{\{a,p\}}(k_r),\textrm{ }q=e^{-\pi\sqrt{r}}, \forall r>0. 
\end{equation} 
Hence finding the expansion of a theta function only requires the knowlege of $Q_{\{a,p\}}(x)$. Instant cases of $Q$ can be found for rational values of $x$ with the routines 'Recognize' or 'RootApproximant' with program Mathematica. An example of evaluation is
$$
\left([1,3;e^{-\pi\sqrt{k_i(1/5)}}]^{*}\right)^6=\frac{1}{90} [-182-\sqrt{689224-148230\cdot 3^{2/3} \sqrt[3]{10}}+
$$
\begin{equation}
+\sqrt{2 \left(92571934 \sqrt{\frac{2}{344612-74115\cdot 3^{2/3} \sqrt[3]{10}}}+74115\cdot 3^{2/3} \sqrt[3]{10}+689224\right)}]
\end{equation}
which is root of the equation
\begin{equation}
45 x^4+364 x^3-21870 x^2-885735=0
\end{equation} 

\[
\]

\centerline{\bf References}\vskip .2in

\noindent

[1]: M. Abramowitz and I.A. Stegun. 'Handbook of Mathematical Functions'. Dover Publications, New York. 1972.

[2]: T. Apostol. Introduction to Analytic Number Theory. Springer Verlang, New York, Berlin, Heidelberg, Tokyo, 1974.

[3]: J.V. Armitage W.F. Eberlein. 'Elliptic Functions'. Cambridge University Press. (2006)

[4]: N.D. Bagis. 'Parametric Evaluations of the Rogers-Ramanujan Continued Fraction'. International Journal of Mathematics and Mathematical Sciences. Vol (2011) 
 
[5]: N.D. Bagis. 'Generalizations of Ramanujan's Continued Fractions'.\\ arXiv:11072393v2 [math.GM] 7 Aug 2012.  

[6]: N.D. Bagis. 'The $w$-modular function and the evaluation of Rogers-Ramanujan continued fraction'. International Journal of Pure and Applied Mathematics. Vol 84, No 1, 2013, 159-169.

[7]: N.D. Bagis. 'Evaluation of the fifth degree elliptic singular moduli'. arXiv:1202.6246v2[math.GM](2012)

[8]: N.D. Bagis. 'On a general sextic equation solved by the Rogers-Ramanujan continued fraction'. arXiv:1111.6023v2[math.GM](2012)

[9]: N.D. Bagis. 'Generalized Elliptic Integrals and Applications'.\\ arXiv:1304.2315v1[math.GM] 4 Apr 2013

[10]: N.D. Bagis and M.L. Glasser. 'Integrals related with Rogers-Ramanujan continued fraction and $q$-products'. arXiv:0904.1641. 10 Apr 2009.

[11]: N.D. Bagis and M.L. Glasser. 'Jacobian Elliptic Functions, Continued Fractions and Ramanujan Quantities'. arXiv:1001.2660v1 [math.GM] 2010.

[12]: N.D. Bagis and M.L. Glasser. 'Conjectures on the evaluation of alternative modular bases and formulas approximating $1/\pi$'. Journal of Number Theory, Elsevier. (2012)

[13]: Bruce C. Berndt. 'Ramanujan`s Notebooks Part III'. Springer Verlag, New York (1991)

[14]: Bruce C. Berndt. 'Ramanujan's Notebooks Part V'. Springer Verlag, New York, Inc. (1998)

[15]: Bruce C. Berndt and Aa Ja Yee. 'Ramanujans Contributions to Eisenstein Series, Especially in his Lost Notebook'. (page stored in the Web). 

[16]: Bruce C. Berndt, Heng Huat Chan, Sen-Shan Huang, Soon-Yi Kang, Jaebum Sohn and Seung Hwan Son. 'The Rogers-Ramanujan Continued Fraction'. J. Comput. Appl. Math., 105 (1999), 9-24.  

[17]: J.M. Borwein and P.B. Borwein. 'Pi and the AGM'. John Wiley and Sons, Inc. New York, Chichester, Brisbane, Toronto, Singapore. (1987)

[18]: D. Broadhurst. 'Solutions by radicals at Singular Values $k_N$ from New Class Invariants for $N\equiv3\;\; mod\;\; 8$'. arXiv:0807.2976 [math-ph], (2008).        

[19]: W. Duke. 'Continued fractions and Modular functions'. Bull. Amer. Math. Soc. (N.S.), 42 (2005), 137-162. 

[20]: I.S. Gradshteyn and I.M. Ryzhik. 'Table of Integrals, Series and Products'. Academic Press (1980).

[21]: Soon-Yi Kang. 'Ramanujan's formulas for the explicit evaluation of the Rogers-Ramanujan continued fraction and theta functions'. Acta Arithmetica. XC.1 (1999)

[22]: E.T. Whittaker and G.N. Watson. 'A course on Modern Analysis'. Cambridge U.P. (1927)

[23]: N.D. Bagis. 'On Generalized Integrals and Ramanujan-Jacobi Special Functions'. arXiv:1309.7247 

\end{document}